\documentclass[11pt]{article}
\usepackage{amsmath,amsthm,amssymb}
\newtheorem{theorem}{Theorem}
\newtheorem{lemma}{Lemma}

\theoremstyle{definition}
\newtheorem{definition}{Definition}
\begin{document}
\title{Bounded operators on the weighted   spaces of
holomorphic functions on the unit Ball in $C^n$}

\author{A. V. Harutyunyan$^1$ and W.Lusky$^2$\\
$^1$Department of Applied Mathematics, Yerevan State University\\ 1 Alex Manookian str., 0025, Yerevan, Armenia\\
anahit@ysu.am\\[2pt]
$^2$Institute of Mathematics, University of Paderborn\\ 100 Warburger str., 33098, Paderborn, Germany\\
lusky@math.uni-paderborn.de}

\maketitle

\begin{abstract}
Assuming that $S$ is the space of functions of regular variation, $\omega\in S$, $0< p<\infty$, a function $f$ holomorphic in $B^n$ is said to be of Besov space $B_p(\omega)$ if
$$\|f\|^p_{B_p(\omega )}=\int_{B^n} (1-|z|^2)^p|Df(z)|^p\frac{\omega(1-|z|)}{(1-|z|^2)^{n+1}}d\nu(z) <+\infty,$$
where $d\nu (z) $ is the volume measure  on $B^n$ and $D $ stands for a fractional derivative  of $f$.

We consider operators on  $B_p(\omega)$ and show, that they are bounded.


\end{abstract}

\section{Introduction and basic constructions}\label{s.1}

 Let  $C^n$  denote the complex Euclidean space of a dimension $n$. For any points $z=(z_1,\ldots,z_n)$, $\zeta=(\zeta_1,\ldots,\zeta_n)$ in $C^n$, we define the inner product as $<z,\zeta>=z_1\overline\zeta_1+\ldots+z_n\overline\zeta_n$ and and note that $|z|^2=|z_1|^2+\ldots+|z_n|^2.$
By $ B^n=\{z\in C^n,\,\,|z|<1\}$ and $C^n:S^n=\{z\in C^n,\,\, |z|=1\}$ we denote the open  unit ball and its boundary, i.e. the unit sphere, in $C^n$. Further, by $H(B^n)$ we denote the set of holomorphic functions on $B^n$ and by $H^\infty (B^n)$ the set of bounded holomorphic functions on $B^n$.

If $f\in H(B^n)$, then  $f(z)=\sum_ma_mz^m$ $(z\in B^n)$, where the sum is taken over all multiindices $m=(m_1,\ldots,m_n)$ with nonnegative integer components $m_k$ and $z^m=z_1^{m_1}\ldots z_n^{m_n}$. Assuming that $|m|=m_1+\ldots+m_n$ and putting $f_k(z)=\sum_{|m|=k}a_mz^m$ for any $k\geq 0$, one can rewrite the Taylor expansion of $f$ as
\begin {equation}\label {e1}
f(z)=\sum_{k=0}^\infty f_k(z),
\end{equation}
which is called homogeneous expansion of $f$, since each $f_k$ is a homogeneous polynomial of the degre $k$.
Further, for a holomorphic function $f$ the fractional differential $D^\alpha $ is defined as
\begin{align*}
D^\alpha f(z)=\sum_{k=0}^\infty (k+1)^\alpha f_k(z),\\
D^\alpha f(\overline{z})=\sum_{k=0}^\infty (k+1)^\alpha f_k(\overline z), \quad k=(k_1,\ldots,k_n),\quad z\in B^n.
\end{align*}
We consider the inverse operator $D^{-\alpha }$ defined in the standard way:
$$D^{-\alpha}D^{\alpha}f(z)=f(z).$$
Particularly, $D^1f(z)=Df(z)$ if $\alpha=1$.

The  following propertes of $D$ are evident

{\em 1}. $DD^\alpha f(z)=D^{\alpha+1}f(z)$

{\em 2}. $D^m(1-<z,\zeta>)^{-\alpha}\preceq (1-<z,\zeta>)^{-\alpha-m}$

By $d\nu$ we denote the volume measure on $B^n$, normalized so that $\nu(B^n)=1$, and by $d\sigma$ the surface measure on $S^n$, normalized so that $\sigma (S^n)=1$. Then following  lemma, the proof of which can be found in \cite{ru} or \cite{zh}, reveals the connection between these measures.
\begin{lemma}\label {l1}
If $f$ is a measureable function with summable modulus over $B^n$, then
$$\int_{B^n}f(z)d\nu(z)=2n\int_0^1r^{2n-1}dr\int_{S^n}f(r\zeta)d\sigma(\zeta).$$
\end{lemma}
\begin{definition}
By $S$ we denote the well-known class of all non-negative measurable  functions $\omega $ on $(0,1)$ with
$$\omega(x)=\exp\biggl\{\int_x^1\frac{\varepsilon (u)}{u}du\biggr\},\,\,x\in (0,1),$$
where $\varepsilon (u)$ is  some measurable, bounded functions on $(0,1)$  and $-\alpha_{\omega}\leq \varepsilon (u)\leq\beta_{\omega}.$
\end{definition}
Note that the functions of $S$ are called {\it{functions of regular variation}} (see \cite{se}). Throughout the paper, we shall assume that  $\omega\in S.$ Throughout the paper the capital Letters $C(...)$ and $C_k$ stand
for different positive constants depending only on the parameters indicated.

We define the holomorphic Besov spaces on the unit ball as follows(see \cite{hl1}).
\begin{definition} Let $\omega \in S,\,\, 0<p<\infty.$
A function $f\in H(B^n)$ is said to be of $B_p(\omega)$ if
$$M^p_f(\omega )=\int_{B^n} (1-|z|^2)^p|Df(z)|^p\frac{\omega(1-|z|)}{(1-|z|^2)^{n+1}}d\nu(z) <+\infty.$$
\end{definition}
We introduce the norm in $H(B^n)$ as $\|f\|_{B_p(\omega)}=M_f(\omega)$ ($|f(0)|$ is not to be added since $Df=0$ implies $f=0$ for a holomorphic function $f$). Besides, it is easy to check that if $p>1$, $n=1$ and $\omega(t)=1$, then $B_p(\omega)$ becomes the classical Besov space (see \cite{ar},\,\,\cite{kk},\,\cite{st}).

\medskip\noindent In particular, for $p=+ \infty $ we shall write $B_\infty (\omega )=B_\omega ,$ where  $B_\omega$ denotes the $\omega-$weighted Bloch space on the ball (see \cite{hl2}).

In \cite{no}, \cite{ls}, \cite{lw}, one can see some other definitions and some characterizations of holomorphic Besov spaces on $B^n.$

Let $1\leq p<\infty$ and let $f\in B_p(\omega)$. Further, let
$m>-n/p-\beta_\omega/p$. Then the function $Df(z)$ has the representation
\begin {equation}\label{e2}
Df(z)=C(\pi,m)\int_{B^n}\frac{(1-|\zeta |^2)^m Df(\zeta)}{(1-<z,\zeta >)^{n+1+m}}d\nu(\zeta),\quad z\in B^n,
\end {equation}
where $C(n,m)=\frac{\Gamma (n+m +1)}{\Gamma(n+1)\Gamma(m+1)}$, follows as a simple consequence of that well known in the one-dimensional case (for details, see \cite{ds}, \cite{zh}).

The following auxiliary lemma will be used.
\begin{lemma}\label{l2}
If $1\geq p<\infty$  and $f\in B_p(\omega)$, then
$$|f(z)|\leq C(\pi,m)\int_{B^n}\frac{(1-|\zeta |^2)^m}{|1-<z,\zeta>|^{n+m}}|D f(\zeta)|d\nu(\zeta)$$
for  $m\in N$ and $m>-n/p-\beta_\omega/p$.
\end{lemma}
{\bf Proof.}
Obviously,  $f(z)=\int_0^1Df(rz)dr$, and by  (\ref {e2}) we get
\begin{align*}
f(z)=&C(\pi,m)\int_0^1\int_{B^n}\frac{(1-|\zeta |^2)^m Df(\zeta)}{(1-r<z,\zeta> )^{n+1+m}}d\nu(\zeta)dr\\
&=C(\pi,m)\int_{B^n}(1-|\zeta |^2)^m Df(\zeta)\int_0^1\frac{dr}{(1-r<z,\zeta> )^{n+1+m}}d\nu(\zeta)\\
&=\widetilde C(\pi,m)\int_{B^n}\frac{(1-|\zeta |^2)^m((1-<z,\zeta>)^{n+\alpha}-1)}{<z,\zeta>(1-<z,\zeta>)^{n+m}}
Df(\zeta)d\nu(\zeta).
\end{align*}
It is clear that $((1-<z,\xi>)^{n+m+1}-1)/<z,\xi>$ is bounded in $B^n$.
Hence the desired statement follows.
 \hfill
 $\square$

\begin{lemma}\label{l3}
Let $\omega\in S$ and let $f\in B_p(\omega)$ for some $0<p\leq 1$. Then
$$
\Biggl(\int_{B^n}|Df(z)|\frac{\omega^{1/p}(1-|z|)}{(1-|z|)^{n}}d\nu(z)\Biggr)^{p}\leq \int_{B^n}|Df(z)|^p\frac{(1-|z|)^p\omega(1-|z|)}{(1-|z|)^{n+1}}d\nu(z)
$$
\end{lemma}
{\bf Proof.}
We have $|Df(z)|=|Df(z)|^p|Df(z)|^{1-p}$. By  Lemma \ref {l2} we get
$$|Df(z)|\leq|Df(z)|^p\frac{\|f\|^{1-p}_{B^p(\omega )}}{\omega^{(1-p)/p}(1-|z|)(1-|z|)^{1-p}}.
$$
Therefore
$$|Df(z)|\frac{(1-|z|)\omega^{1/p}(1-|z|)}{(1-|z|)^{n+1}}\leq |Df(z)|^p\|f\|^{1-p}_{B^p(\omega )}\frac{\omega(1-|z|)(1-|z|)^{p}}{(1-|z|)^{n+1}},
$$
and by integration over $B^n$ we get
$$\int_{B^n}|Df(z)|\frac{(1-|z|)\omega^{1/p}(1-|z|)}{(1-|z|)^{n+1}}d\nu(z)\leq \|f\|^{1-p}_{B^p(\omega )}\int_{B^n}|Df(z)|^p\frac{\omega(1-|z|)(1-|z|)^{p}}{(1-|z|)^{n+1}},
$$

The proof is completed.
 \hfill
 $\square$

 \begin{lemma}\label{l4}
 Let  $\omega\in S ,\,\,\,\alpha+1-\beta_\omega >0,$ and
 $\beta-\alpha>\alpha_\omega$. Then
 $$\int_{B^n}\frac{(1-|\zeta|^2)^\alpha\omega(1-|\zeta|)}{|1-<z, w>|^{\beta +n+1}}d\nu (\zeta)
 \leq C(\alpha,\beta,\omega)\frac{\omega (1-|z|^2)}{(1-|z|^2)^{\beta-\alpha}}.
 $$
 \end{lemma}

 {\bf Proof.}
By Lemma \ref{l1} for $\beta>0$ we get
 $$
 \int_{B^n}\frac{(1-|\zeta|^2)^\alpha\omega(1-|\zeta|)}{|1-<z,\zeta>|^{\beta+n+1}}d\nu(\zeta)=2n\int_0^1r^{2n-1}(1-r^2)^\alpha\omega(1-r)dr\times$$
$$
\int_{S^n}\frac{d\sigma(\zeta)}{|1-<z,\zeta>|^{\beta+n+1}}\leq
2n\int_0^1r^{2n-1}\frac{(1-r^2)^\alpha\omega(1-r)}{(1-r|z|)^{\beta+1}}dr.$$

In the last inequality we have used Theorem 1. 12 from \cite{zh}.

The problem is to estimate the last one dimensional integral. To this end we have

$$\int_0^1\frac{(1-r^2)^{\alpha}\omega (1-r) dr}{(1-r^|z|)^{\beta+1}}\leq\int_0^1\frac{u^\alpha\omega (u) du}{(1-|z|+u|z|)^{\beta+1}}$$
$$=\left\lbrace \int_0^{1-|z|}\frac{u^\alpha\omega (u)
    du}{(1-|z|+u|z|)^{\beta+1}}+\int_{1-|z|}^1\frac{u^\alpha\omega
    (u) du}{(1-|z|+u|z|)^{\beta+1}}\right\rbrace =I_1+I_2.$$

First we  estimate the integral $I_1.$
$$I_1\leq \int_0^{1-|z|}\frac {u^\alpha\omega (u)}{(1-|z|)^{\beta+1}}du=\frac{1}{(1-|z|)^{\beta+1}}\int_0^{1-|z|}u^\alpha\omega (u)du=$$
$$\frac {(\alpha+1)^{-1}}{(1-|z|)^{\beta+1}}
\left[ (1-|z|)^{\alpha+1}\omega (1-|z|)+\int_0^{1-|z|}u^\alpha\omega (u)\varepsilon (u)du\right]$$
As a result we get
$$(\alpha+1)\int_0^{1-|z|}u^\alpha\omega (u)du= (1-|z|)^{\alpha+1}\omega (1-|z|)+\int_0^{1-|z|}u^\alpha\omega (u)\varepsilon (u)du, $$
and
$$\int_0^{1-|z|}(\alpha+1-\varepsilon (u))u^\alpha\omega (u)du= (1-|z|)^{\alpha+1}\omega (1-|z|).$$
On  the other hand
$$\alpha+1-\beta_\omega\leq \alpha+1-\varepsilon (u)$$
which yields
$$(\alpha+1-\beta_\omega )\int_0^{1-|z|}u^\alpha\omega (u)du\leq (1-|z<)^{\alpha+1}\omega(1-|z|)$$
or
\begin{equation}\label{e3}I_1\leq C(\alpha,\beta,\omega)\frac{\omega
    (1-|z|)}{(1-|z|)^{\beta-\alpha}}\end{equation}

Now we want to  estimate  $I_2$. For  $|z|\geq 1/2$ we have
$$I_2\leq 2^{\beta-1}\int_{1-|z|}^1\frac{\omega (u)}{u^{\beta-\alpha+1}}=\frac{2^{\beta}}{\alpha-\beta-1}\int_{1-|z|}^1\omega (u)du^{\alpha-\beta-1}$$
$$=\frac{2^{\beta-1}}{\beta-\alpha+1}\left[ \frac{\omega (1-|z|)}{(1-|z|)^{\beta-\alpha-1}}+\omega (1)+\int_{1-|z|}^1\frac{\omega (u)\varepsilon (u)}{u^{\beta-\alpha}}du\right] $$
Then
$$\int_{1-|z|}^1\frac{\omega (u)}{u^{\beta-\alpha}}du= \frac{\omega (1-|z|)}{(1-|z|)^{\beta-\alpha-1}}-\omega (1)-\int_{1-|z|}^1\frac{\omega (u)\varepsilon (u)}{u^{\beta-\alpha}}du$$
and it  follows that
$$\int_{1-|z|}^1\left( 1+\frac{\varepsilon (u)}{\beta-\alpha-1}\right) \frac{\omega (u)}{u^{\beta-\alpha}}du= \frac{\omega (1-|z|)}{(1-|z|)^{\beta-\alpha-1}}-\omega (1)\leq\frac{\omega (1-|z|)}{(1-|z|)^{\beta-\alpha-1}}.$$
Then the inequality  $$1+\frac{\varepsilon (u)}{\beta-\alpha-1}\geq 1-\frac{\alpha_\omega}{\beta-\alpha-1} >0,$$ gives us
\begin{equation}\label{e4}\int_{1-|z|}^1\frac{\omega (u)}{u^{\beta-\alpha}}du\leq C(\alpha,\beta,\omega)\frac{\omega (1-|z|)}{(1-|z|)^{\beta-\alpha-1}}.\end{equation}
Summing up, from (\ref{e3})  and (\ref{e4}) we get the proof of Lemma \ref{l4}.

\hfill
$\square$

\begin{lemma}\label{l5} The following statement is true
  $$D(fg)= gDf+ fDg-fg$$.
 \end{lemma}
{\bf Proof.}   To this end we define the radial derivative $R$ of $f$ as follows
\[ (Rf)(z) = \sum_{k=1}^{ \infty} k f_k(z), \ \ \ \ z \in B^n \]
or equivalently
\[ (Rf)(z)= \sum_{k=1}^n z_k \frac{ \partial f(z)}{ \partial z_k}, \ \ \  z \in B^n. \]
It is easy to note that $R(fg) = gR(f) + f R(g)$ and for the operator $D$ we have
$Df = f + Rf$. Combining the last equalities we get
\[ D(fg) = fg + g R(f) + f R(g). \]
On the other hand we have $R(f) = Df -f$ and $R(g) = Dg-g$. Hence $D(fg) = gDf + fDg -fg$.

 \hfill
 $\square$

\section{Bounded Operators on $B_p(\omega)$}\label{s.2}

In this Section first we consider the following operator
 \begin{eqnarray*}
T_{\bar{h}}^{ \alpha}(f)(z) & =
 & \int_{B^n} \frac{(1-|\xi|^2)^{ \alpha} \overline{h(\xi)}f(\xi)} {(1- \langle z, \xi \rangle)^{n+ \alpha + 1}}d \vartheta(\xi),\,\,\alpha>-1.
\end{eqnarray*}

\begin{theorem} Let $0 < p < \infty$, $ h \in H^1(B^n)$. Then 

1. if $T_{\bar{h}}^{ \alpha}$ is bounded on $B_p( \omega)$ then  $h \in H^{ \infty}( B^n)$.

2. conversely,

 a) if $1\leq p<\infty$ and $h\in H^\infty(B^n)$ then $T_{\bar{h}}^{ \alpha}: B_p(\omega)\rightarrow B_p(\omega)$

b) if $0<p<1$ and $h\in H^\infty(B^n)$ then $T_{\bar{h}}^{ \alpha}: B_p(\omega)\rightarrow B_p(\omega^*), $ where $\omega^*(t)=t^{(\alpha+m+1)(1-p)}\omega(t)$ and $m>-n/p-\beta_\omega/p$.
\end{theorem}
{\bf Proof.}
1. Let $T_{\bar{h}}^{ \alpha}$ be bounded on $B_p( \omega)$. We take
\[ f_{\tau}(\xi) = \frac{1}{(1- \langle \xi, \tau \rangle )^{ \alpha+n+1}}, \]
where $ \tau \in [0,1]^n$ is a parameter. We calculate
\begin{eqnarray*}
T_{\bar{h}}^{ \alpha}(f)(z) & = & \int_{B^n} \frac{(1-|\xi|^2)^{ \alpha} \overline{h(\xi)} d \vartheta(\xi)}{(1- \langle z, \xi \rangle)^{n+ \alpha + 1}(1- \langle \xi, \tau \rangle^{n+ \alpha+1}} \\
                            & = & \overline{ \int_{B^n} \frac{(1-|\xi|^2)^{ \alpha} h(\xi) d \vartheta(\xi)}{(1- \langle  \xi,z \rangle)^{n+ \alpha + 1}(1- \langle \tau,\xi  \langle)^{n+ \alpha+1}}} 
                            =  \frac{ \overline{h(z)}}{(1-\langle z, \tau \rangle)^{n+ \alpha + 1}}
\end{eqnarray*}
Then we get

\[ ||T_{\bar{h}}(f_{\tau}||_{B_p( \omega)} = | h(z) \cdot ||f_{\tau}||_{B_p( \omega)} \leq ||T_{\bar{h}}|| \cdot ||f||_{B_p( \omega)} \]
and hence $|h( \tau)| \leq ||T_{\bar{h}}||$.
Replacing $f_{\tau}(z)$ by $f_{\tau}(e^{i \theta}z)$ we get
$|h( \tau e^{i \theta})| \leq ||T_{\bar{h}}||$ which implies $h \in H^{ \infty }(B^n)$.

2. Conversely,
a) let $p\geq 1$ and  $h \in H^{ \infty}(B^n)$. We show that $T_{\bar{h}}^{ \alpha}(f) \in
B_p( \omega)$ for any $f \in B_p( \omega)$. To this end by Lemma \ref{l2} we use the inequality
\[ |f(\xi) | \leq C(\pi,m)\int_{B^n} \frac{(1- |t|^2)^m |Df(t)|}{|1-  \xi,t \rangle|^{m+n}} d \vartheta(t) \]
which implies that
\[ |f(\xi) |^p \leq \frac{C(\pi,m)}{(1-|\xi|^2)^{(m-1)p/q}} \int_{B^n} \frac{(1- |t|^2)^{mp} |Df(t)|^p}{|1- \langle \xi, t \rangle|^{m+n}} d \vartheta(t) \]
Then  for $ p > 1$, by Holders inequality, we get,
\begin{eqnarray*}
|T_{\bar{h}}^{ \alpha}f(z)|^p & \leq & C(\pi,m)\left( \int_{B^n} \frac{(1-|\xi|^2)^{\alpha}|f(\xi)| \cdot| \overline{h(\xi)|} d \vartheta(\xi)}{|1- \langle z, \xi \rangle|^{n+2+ \alpha}} \right)^p \\
                              & \leq & C(\pi,m)\frac{||h||_{ \infty}}{(1-|z|^2)^{p/q}}
\int_{B^n}  \frac{(1-|\xi|^2)^{\alpha}|f(\xi)|^p  d \vartheta(\xi)}{|1- \langle z, \xi \rangle|^{n+2+ \alpha}} .
\end{eqnarray*} Then we have
\begin{eqnarray*}
I\equiv C(\pi,m)\int_{B^n} |DT_{\bar{h}}^{ \alpha}f(z)|^p \frac{ \omega(1-|z|) d \vartheta(z)}{(1-|z|^2)^{n+1-p}} \\
  \leq C(\pi,m) \int_{B^n} (1-|t|^2)^{mp}|Df(t)|^p \int_{B^n} \frac{(1-|\xi|^2)^{\alpha-(m-1)p/q}}{|1-\langle \xi,t \rangle|^{m+n} } \\ \int_{B^n} \frac{ \omega(1-|z|)d \vartheta(z) d\vartheta(\xi) d\vartheta(t)}{(1-|z|^2)^{n+1-p+p/q}|1- \langle z,\xi \rangle|^{n+2+ \alpha}}
  \end{eqnarray*}
 Using Lemma \ref{l4} we obtain furthermore
 \begin{eqnarray*}
 & &I \leq C(\pi,m)\int_{B^n} (1-|t|^2)^{mp}|Df(t)|^p\int_{B^n} \frac{ (1-|\xi|^2)^{ \alpha}\omega(1-|\xi|)d \vartheta(z) d\vartheta(\xi) }{|1-|\xi|^2)^{n+1+ \alpha +(m-1)p/q}|1- \langle \xi,t \rangle|^{n+m}}  \\
  & & \leq C(\pi,m) ||h||_{ \infty}\int_{B^n} \frac{ (1-|t|^2)^{ mp}|Df(t)|^p\omega(1-|t|)d\vartheta(t) }{(1-|t|^2)^{(m-1)p/q+m+n}} \\& & = C(\pi,m) ||h||_{ \infty}\int_{B^n}\frac{|Df(\xi)|^p\omega(1-|\xi|)d\vartheta(\xi) }{|1-|\xi|^2)^{n+1-p}} 
 \leq C(\pi,m)||f||_{B_p( \omega)}^p ||h||_{ \infty}.
\end{eqnarray*}

Let now $p=1$. We have
\begin{eqnarray*}
\lefteqn{\int_{B^n} |DT_{\bar{h}}^{ \alpha}f(z)| \frac{ \omega(1-|z|) d \vartheta(z)}{(1-|z|^2)^{n}} } \\
 & & \leq C(\pi,m)||h||_{\infty} \int_{B^n} \frac{(1-|w|^2)^{m}|Df(w)|}{|1- \langle \xi, w \rangle |^{m+n}} \int_{B^n} \frac{(1-|\xi|^2)^{\alpha}}{|1-\langle z,w \rangle|^{n+ \alpha+2}} \\
 & & \hspace{6cm}  \int_{B^n} \frac{ \omega(1-|\xi|)d \vartheta(z) d\vartheta(\xi) d\vartheta(w)}{|1-|\xi|^2)^{n}} \\
 & & = C(\pi,m)||h||_{ \infty} \int_{B^n} (1-|w|^2)^{m}|Df(w)|\int_{B^n} \frac{ (1-|\xi|^2)^{ \alpha} }{|1- \langle \xi,w \rangle|^{n+m}} \\
 & &  \hspace{6cm} \int_{B^n}  \frac{ \omega(1-|z|)d \vartheta(z) d\vartheta(\xi) d\vartheta(w)}{(1-|z|^2)^{n}|1- \langle z, \xi \rangle |^{n+2+ \alpha}} \\
 & & C(\pi,m)\leq ||h||_{ \infty}\int_{B^n}  (1-|w|^2)^{ m}|Df(w)| \int_{B^n} \frac{(1- |\xi|^2)^{ \alpha}\omega(1-|\xi|)d\vartheta(\xi) d \vartheta(w)}{|1- \langle \xi, w \rangle|^{m+n}(1-|\xi|^2)^{n+1+ \alpha}} \\
 & & C(\pi,m)\leq   ||h||_{ \infty}\int_{B^n}  (1-|w|^2)^{ m}|Df(w)| \frac{ \omega(1-|w|)(1-|w|^2)^{ \alpha} \vartheta(w) }{|1-|w|^2)^{n+m+ \alpha}} \\
 & & = \int_{B^n} \frac{|Df(w)| \omega(1-|w|)}{(1-|w|^2)^n} d \vartheta(w) =C(\pi,m) ||f||_{B_p( \omega)} ||h||_{ \infty}.
\end{eqnarray*}

b) Let $0<p<1$.
Using Lemma \ref{l3} we get
\[ |f(\xi) |^p \leq C(\pi,m) \int_{B^n} \frac{(1- |t|^2)^{mp+p} |Df(t)|^p}{|1- \langle \xi, t \rangle|^{m+n+1}} d \vartheta(t) \]
and
\[|DT_{\bar{h}}^{ \alpha}f(z)|^p  \leq C(\pi,m) ||h||^p_{ \infty}\int_{B^n}  \frac{(1-|\xi|^2)^{p(\alpha+1)}|f(\xi)|^p  d \vartheta(\xi)}{|1- \langle z, \xi \rangle|^{n+2+ \alpha}}. \]
As in the case of $p>1$ by Lemma \ref{l4} we  get
\begin{eqnarray*}
\lefteqn{\int_{B^n} |DT_{\bar{h}}^{ \alpha}f(z)|^p \frac{ \omega^*(1-|z|) d \vartheta(z)}{(1-|z|^2)^{n+1-p}} } \\
 & & C(\pi,m)\leq \int_{B^n} (1-|t|^2)^{mp+p}|Df(t)|^p \int_{B^n} \frac{(1-|\xi|^2)^{p(\alpha+1)}}{|1-\langle \xi,t \rangle|^{m+n+1} } \\
 & & \hspace{4.5cm} \int_{B^n} \frac{ \omega^*(1-|z|)d \vartheta(z) d\vartheta(\xi) d\vartheta(t)}{(1-|z|^2)^{n+1-p}|1- \langle z,\xi \rangle|^{n+2+ \alpha}} \\
 & & C(\pi,m)\leq \int_{B^n} (1-|t|^2)^{mp+p}|Df(t)|^p\int_{B^n} \frac{ (1-|\xi|^2)^{p( \alpha+1)}\omega^*(1-|\xi|)d \vartheta(z) d\vartheta(\xi) }{(1-|\xi|^2)^{n+1-p+\alpha+2}(1- \langle \xi,t \rangle|^{n+m+1}}  \\
 & & C(\pi,m)\leq ||h||_{ \infty}\int_{B^n} \frac{ (1-|t|^2)^{ mp+p}|Df(t)|^p\omega^*(1-|t|)d\vartheta(t) }{(1-|t|^2)^{m+n+2-2p-p\alpha+\alpha}} \\
 & & =  ||h||_{ \infty}\int_{B^n} \frac{ |Df(\xi)|^p\omega(1-|\xi|)d\vartheta(\xi) }{|1-|\xi|^2)^{n+1-p}} 
  \leq C(\pi,m)||f||_{B_p( \omega)}^p ||h||_{ \infty}.
\end{eqnarray*}

\hfill
 $\square$
The next theorem is about boundedness $M_h$
\begin{theorem}
Let $H^{ \infty}(B^n)$. Then $M_h$ is a bounded operator $B_p( \omega) \rightarrow B_p( \omega)$.
\end{theorem}

{\bf Proof.} Using Lemma \ref{l5} we show that
\[ \int_{B^n} |Df(z)|^p |g(z)|^p \frac{ \omega(1-|z|)}{(1-|z|^2)^{n+1-p}}d \vartheta(z) < \infty \hspace{2cm} (3) \]
\[ \int_{B^n} |f(z)|^p |Dg(z)|^p \frac{ \omega(1-|z|)}{(1-|z|^2)^{n+1-p}}d \vartheta(z) < \infty \hspace{2cm} (4) \]
\[ \int_{B^n} |f(z)|^p |g(z)|^p \frac{ \omega(1-|z|)}{(1-|z|^2)^{n+1-p}}d \vartheta(z) < \infty \hspace{2.5cm} (5) \]
The proof of (3) is evident.

Proof of (4). First we show that
\[ |Dg(z) | \leq \frac{C\|g\|_{\infty}}{1-|z|^2}, \ \ \ z \in B^n, \]
if $g \in H^{ \infty}(B^n)$. To this end we take the ball $B^n(z) = \{ w \in B^n, |w-z| < n - |z|/2 \}$ and use the Cauchy inequality. In the case $p > 1$ we have
\begin{eqnarray*}
\lefteqn{ |f(z)|^p \leq \left(C(\pi,m) \int_{B^n} \frac{(1-|\xi|^2)^m|Df(\xi)| d\vartheta(\xi)}{|1-\langle m, \xi \rangle|^{m+n}} \right)^p} \\
 & & \leq \frac{C(\pi,m)}{(1-|z|^2)^{ \gamma p/q}} \int_{B^n} \frac{(1-|\xi|^2)^{m-1- \gamma}(1-|\xi|^2)^{p+ \gamma p}|Df(\xi)|^p d \vartheta(\xi)}{|1-\langle z, \xi\rangle |^{m+n}}
\end{eqnarray*}
Then for (4) we get
\begin{eqnarray*}
\lefteqn{ \int_{B^n} |f(z)|^p |Dg(z)|^p \frac{ \omega(1-|z|)}{(1-|z|^2)^{n+1-p}}d \vartheta(z) } \\
 & & \leq C(\pi,m)\|g\|^p_{\infty} \int_{B^n} (1- | \xi|^2)^{m+(p-1)( \gamma +1)}|Df(\xi)|^p \int_{B^n}
\frac{ \omega(1-|z|^2) d \vartheta(z) d \vartheta(\xi)}{|1- \langle z, \xi \rangle|^{m+n}(1-|z|^2)^{n+2-p+ \gamma p/q}} \\
 & & \leq C(\pi,m)\|g\|_{\infty}\int_{B^n} (1- | \xi|^2)^{m+(p-1)( \gamma +1)} \frac{|Df(\xi)|^p \omega(1- | \xi|)}{(1- |\xi|)^{m+n-p+ \gamma p/q}} d \vartheta(\xi) \\
 & & = C(\pi,m)\|g\|^p_{\infty}\int_{B^n} \frac{|Df( \xi)|^p \omega(1- |\xi|) d \vartheta(\xi) }{(1-|\xi|)^{n+1+-2p}} \\
 & & = \int_{B^n} \frac{|Df( \xi)|^p \omega(1- |\xi|)(1-|\xi|)^{p} d \vartheta(\xi) }{(1-|\xi|^2)^{n+1-p}}   \leq ||f||_{B_p( \omega)}^pC(\pi,m)\|g\|^p_{\infty} \\
\end{eqnarray*}
Proof of (5).
\begin{eqnarray*}
\lefteqn{ \int_{B^n} |f(z)|^p |g(z)|^p \frac{ \omega(1-|z|)}{(1-|z|^2)^{n+1-p}}d \vartheta(z) } \\
 & &C(\pi,m)\|g\|^p_{\infty} \leq  \int_{B^n} (1- | \xi|^2)^{m+(p-1)( \gamma +1)}|Df(\xi)|^p \int_{B^n}
\frac{ \omega(1-|z|) d \vartheta(z) d \vartheta(\xi)}{|1- \langle z, \xi \rangle |^{m+n}(1-|z|^2)^{n+1-p+ \gamma p/q}} \\
 & & \leq \int_{B^n} (1- | \xi|^2)^{m+(p-1)( \gamma +1)} \frac{|Df(\xi)|^p \omega(1- | \xi|)}{(1- |\xi|)^{n+1-p+ \gamma p/q+m-1}} d \vartheta(\xi) \\
& & = C(\pi,m)\|g\|^p_{\infty}\int_{B^n} \frac{|Df( \xi)|^p \omega(1- |\xi|)(1-|\xi|)^{p} d \vartheta(\xi) }{(1-|\xi|^2)^{n+1-p}}   \leq C(\pi,m)\|g\|^p_{\infty} ||f||^p_{B_p( \omega)} 
\end{eqnarray*}
Let $0 < p \leq 1$.
Then by Lemma 4
\[ |f(z)|^p \leq C(\pi,m) \int_{B^n} \frac{(1-|\xi|^2)^{mp+(n+1)p}|Df(\xi)|^p d\vartheta(\xi)}{|1-\langle m, \xi \rangle|^{(m+1)p}(1-|\xi|^2)^{n+1}} \]
Then for (4) we get
\begin{eqnarray*}
\lefteqn{
\int_{B^n} |f(z)|^p |Dg(z)|^p \frac{ \omega(1-|z|)}{(1-|z|^2)^{n+1-p}}d \vartheta(z) } \\
 & &  \leq C(\pi,m)\|g\|^p_{\infty} \int_{B^n} \frac{(1- | \xi|^2)^{mp+(n+1)p}|Df(\xi)|^p}{(1-|\xi|^2)^{n+1}} \int_{B^n}
\frac{ \omega(1-|z|) d \vartheta(z) d \vartheta(\xi)}{|1- \langle z, \xi \rangle |^{(m+1)p}(1-|z|^2)} \\
 & & \leq C(\pi,m)\|g\|^p_{\infty}\int_{B^n} \frac{ (1- | \xi|^2)^{mp+(n+1)p} |Df(\xi)|^p \omega(1- | \xi|)}{(1- |\xi|^2)^{n+1}(1-|\xi|^2)^{(m+1)p}} d \vartheta(\xi) \\
 & & = C(\pi,m)\|g\|^p_{\infty}\int_{B^n} \frac{|Df( \xi)|^p \omega(1- |\xi|)(1-|\xi|^2)^{n(p+1)-p} d \vartheta(\xi) }{(1-|\xi|^2)^{n+1-p}}=C(\pi,m)\|g\|^p_{\infty}\|f\|_{B_p(\omega)}
\end{eqnarray*}
Finally we obtain (5)
\begin{eqnarray*}
 \lefteqn{\int_{B^n} |f(z)|^p |g(z)|^p \frac{ \omega(1-|z|)}{(1-|z|^2)^{n+1-p}}d \vartheta(z)} \\
 & & \leq C(\pi,m)\|g\|^p_{\infty}\int_{B^n} \frac{(1-|\xi|^2)^{mp+(n+1)p}|Df(\xi)|^p}{(1-|\xi|^2)^{n+1}} \frac{ \omega(1-|\xi|) d \vartheta(\xi)}{(1-|\xi|^2)^{(m+1)p-n-1-p}} \\
 & & \leq C(\pi,m)||f||_{B_p( \omega)}^p ||g||_{ \infty}^p
\end{eqnarray*}

Summing ab, we get the proof of theorem.

\hfill
 $\square$

\end{document}